\documentclass[11pt,a4paper]{article}

\usepackage{theorem,enumerate}
\usepackage{amsmath,latexsym,amssymb,amsfonts}
\usepackage{eucal}
\usepackage{mathrsfs}
\usepackage{array}
 \usepackage{graphicx, color}

\newcounter{Scounter}
\theorembodyfont{\normalfont\slshape}
\newtheorem{thm}{Theorem}

\newtheorem{cor}[thm]{Corollary}
\newtheorem{prop}[thm]{Proposition}
\newtheorem{lemma}[thm]{Lemma}
\newtheorem{definition}[thm]{Definition}

\newtheorem{claim}{Claim}

\newcommand{\proof}{\medbreak\noindent\textit{Proof.}\quad}

\newcommand{\qed}{{$\quad\square$\vs{3.6}}}

\newcommand{\vs}[1]{\vspace*{#1 mm}}

\def\M{{ \mathcal{M}}}

\bfseries\normalfont

\makeatletter
\def\thanks#1{%
   \footnotemark
   \edef\@tempa{\noexpand\noexpand\noexpand\footnotetext[\the\c@footnote]}%
   \toks@\expandafter{\@thanks}%
   \toks\tw@{{#1}}
   \xdef\@thanks{\the\toks@\@tempa\the\toks\tw@}}
\makeatother

\begin{document}

\title{A Sufficient condition for DP-4-colorability}

\author{
Seog-Jin Kim\thanks{Department of Mathematics Education, Konkuk University,
Korea.
e-mail: {\tt skim12@konkuk.ac.kr
}}\thanks{
This research was supported by Basic Science Research Program through the National Research Foundation of Korea(NRF) funded by the Ministry of Education(NRF-2015R1D1A1A01057008).}
and 
Kenta Ozeki\thanks{Faculty of Environment and Information Sciences,
Yokohama National University,
Japan. 
e-mail: {\tt ozeki-kenta-xr@ynu.ac.jp}
}\thanks{
This work was supported by JST ERATO Grant Number JPMJER1201, Japan.
}
}

\date{August 22, 2017}
\maketitle

\begin{abstract}
\emph{DP-coloring} of a simple graph is
a generalization of list coloring, and also a generalization of signed coloring 
of signed graphs.  It is known that for each $k \in \{3, 4, 5, 6\}$, every planar graph without $C_k$ is 4-choosable.  Furthermore, Jin, Kang, and Steffen \cite{JKS} showed that for each $k \in \{3, 4, 5, 6\}$, every signed planar graph without $C_k$ is signed 4-choosable. In this paper, we show that for each $k \in \{3, 4, 5, 6\}$, every planar graph without $C_k$ is 4-DP-colorable, which is an extension of the above results.
\end{abstract}


\section{Introduction}

\subsection{List-coloring}

We denote by $[k]$ the set of integers from $1$ to $k$.
A $k$-coloring of a graph $G$
is a mapping $f : V(G) \rightarrow [k]$
such that $f(u) \not= f(v)$ for any $uv \in E(G)$.
(In this paper,
we always use the term ``coloring''
as a \emph{proper} coloring.)
The minimum integer $k$
such that $G$ admits a $k$-coloring 
is called the \emph{chromatic number of $G$},
and denoted by $\chi(G)$.

A \emph{list assignment} $L : V(G) \rightarrow 2^{[k]}$ of $G$
is a mapping that assigns a set of colors to each vertex.
A coloring $f: V(G) \rightarrow Y$ where $Y$ is a set of colors
is called an \emph{$L$-coloring} of $G$ if
$f(u) \in L(u)$ for any $u \in V(G)$.
A list assignment $L$ is called a \emph{$t$-list assignment} if 
$|L(u)| \geq t$ for any $u \in V(G)$.
A graph $G$ is said to be \emph{$t$-choosable}
if $G$ admits an $L$-coloring for each $t$-list assignment $L$,
and 
the \emph{list-chromatic number} or the \emph{choice number} of $G$,
denoted by $\chi_{\ell}(G)$,
is the minimum integer $t$
such that $G$ is $t$-choosable.


Since a $k$-coloring corresponds to an $L$-coloring
with $L(u) = [k]$ for any $u \in V(G)$,
we have $\chi(G) \leq \chi_{\ell}(G)$.
It is well-known that 
there are infinitely many graphs $G$ 
satisfying $\chi (G) < \chi_{\ell}(G)$,
and the gap can be  arbitrarily large.

Thomassen \cite{Thomassen} showed that every planar graph is 5-choosable, and Voigt \cite{Voigt} showed that there are planar graphs which are not 4-choosable. 
Thus finding sufficient conditions for planar graphs to be 4-choosable is an interesting problem.

Let $C_k$ be the cycle of length $k$.
Lam, Xu, and Liu \cite{Lam} showed that every planar graph without $C_4$ is 4-choosable.  And for each $k \in \{3, 5, 6\}$, it is known that if $G$ is a planar graph without $C_k$, then $G$ is 3-degenerate \cite{FJMS, LW02} (see Lemma \ref{356-lemma} in Section \ref{main-proof-sec}),
where a graph is said to be \emph{$\ell$-degenerate}
if its arbitrary subgraph contains a vertex of degree at most $\ell$.  Thus we can conclude that for each 
$k \in \{3, 5, 6\}$, if $G$ is a planar graph without $C_k$, then $\chi_{\ell} (G) \leq 4$.  Thus with the result in \cite{Lam}, we have the following.

\begin{thm} \label{list-noC4}
For each 
$k \in \{3, 4,  5, 6\}$, if $G$ is a planar graph without $C_k$, then $\chi_{\ell} (G) \leq 4$.
\end{thm}

Notice that 
for $k \in \{3, 4,  5, 6\}$,
Theorem \ref{list-noC4} is best possible.
That is,
for $k \in \{3, 4,  5, 6\}$,
there is a planar graph $G$ without $C_k$
but with $\chi_{\ell} (G) \geq 4$.
See \cite{Voigt} for $k =3$ and \cite{Voigt2} for $k = 4,5$.
For $k = 6$, complete graph $K_4$ has no $C_6$ but $\chi_{\ell}(K_4) = 4$.

\subsection{Singed colorings of signed graphs}

A \emph{signed graph} $(G,\sigma)$
is a pair of a graph $G$ 
and a mapping $\sigma: E(G) \rightarrow \{1, -1\}$,
which is called a \emph{sign}.
For an integer $k$,
let 
$$N_k = 
\begin{cases}
\{0, \pm 1, \dots , \pm r\} & \text{if $k$ is an odd integer with $k = 2r +1$,} \\
\{\pm 1, \dots , \pm r\} & \text{if $k$ is an even integer with $k = 2r$.}
\end{cases}
$$
Note that $|N_k| = k$.
A \emph{signed $k$-coloring} of a signed graph $(G,\sigma)$
is a mapping $f : V(G) \rightarrow N_k$
such that $f(u) \not= \sigma(uv) f(v)$
for each $uv \in E(G)$.
The minimum integer $k$
such that a signed graph $(G,\sigma)$ admits a signed $k$-coloring
is called the \emph{signed chromatic number of $(G,\sigma)$}.
This was first defined by Zaslavsky \cite{Zaslavsky}
with slightly different form,
and then modified by M\'{a}\v{c}ajov\'{a}, Raspaud, and \v{S}koviera \cite{MRS}
to the above form
so that it would be a natural extension of an ordinary vertex
coloring.

Now we define signed list-colorings of signed graphs.  Given a signed graph $(G, \sigma)$, a list-assignment of $(G, \sigma)$ is a function $L$ defined on $V(G)$ such that $L(v) \subset \mathbb{Z}$ for each $v \in V(G)$. 
A \emph{signed $L$-coloring} $\phi$ of $(G, \sigma)$ is a signed coloring such that $\phi(v) \in L(v)$ for each $v \in V(G)$. 
A signed graph $(G, \sigma)$ is called \emph{signed $k$-choosable} if it admits an $L$-coloring for every 
list-assignment $L$ with
$L(v) \subset \mathbb{Z}$ and $|L(v)| \geq k$ for each $v \in V(G)$. 
 The {\em signed choice number} of $(G, \sigma)$ is the minimum number $k$ such that $(G, \sigma)$ is signed $k$-choosable.

Note that if $\sigma (uv) = +1$ for all edge $uv$ in $E(G)$, then the signed choice number of $(G, \sigma)$ is the same as the choice number of the graph $G$.
Jin, Kang, and Steffen \cite{JKS} showed the following theorem.

\begin{thm} \cite{JKS} \label{signed-noC4}
For each $k \in \{3, 4, 5, 6\}$, every signed planar graph without $C_k$ is signed $4$-choosable.
\end{thm}

Note that Theorem \ref{signed-noC4} is an extension of Theorem \ref{list-noC4}.
We will give a further extension of Theorem \ref{signed-noC4},
together with a proof different from that in \cite{JKS}.


\subsection{DP-coloring}
In order to consider some problems on list chromatic number,
Dvo\v{r}\'{a}k and Postle \cite{DP} considered 
a generalization of a list-coloring.
They call it a \emph{correspondence coloring},
but we call it a \emph{DP-coloring},
following Bernshteyn, Kostochka and Pron \cite{BKP}.

Let $G$ be a  graph 
and $L$ be a list assignment of $G$.
For each edge $uv$ in $G$,
let $M_{L,uv}$ be a matching
between $\{u\} \times L(u)$ and $\{v\} \times L(v)$.
With abuse of notation,
we sometimes regard $M_{L,uv}$ 
as a bipartite graph
between $\{u\} \times L(u)$ and $\{v\} \times L(v)$
of maximum degree at most 1.

\begin{definition}
Let $\M_{L} = \big\{M_{L,uv} : uv \in E(G) \big\}$,
which is called a \emph{matching assignment over $L$}.
Then a graph $H$ is said to be the \emph{$\M_{L}$-cover} of $G$
if it satisfies all the following conditions:
\begin{enumerate}[{\upshape (i)}]
\item
The vertex set of $H$ is 
$\bigcup_{u \in V(G)} \big(\{u\} \times L(u)\big) 
= \big\{(u,c): u \in V(G), \ c \in L(u)\big\}$.

\item
For any $u \in V(G)$,
the set $\{u\} \times L(u)$ induces a clique in $H$.
\item
For any edge $uv$ in $G$,
$\{u\} \times L(u)$
and 
$\{v\} \times L(v)$
induce in $H$ the graph obtained from $M_{L,uv}$
by adding those edges defined in (ii).

\end{enumerate}
\end{definition}

\begin{definition}
An $\M_{L}$-coloring of $G$ is an independent set $I$ in 
the $\M_{L}$-cover with $ |I| = |V(G)|$.
The \textit{DP-chromatic number},
denoted by $\chi_{\text{DP}}(G)$,
is the minimum integer $t$
such that $G$ admits an $\M_{L}$-coloring
for each $t$-list assignment $L$ and each matching assignment $\M_{L}$ over $L$.  We say that a graph $G$ is \emph{DP-$k$-colorable} if $\chi_{DP}(G) \leq k$.
\end{definition}

Note that 
when $G$ is a simple graph and 
$$M_{L,uv} = \big\{(u,c)(v,c): c \in L(u) \cap L(v)
\big\}
$$
for any edge $uv$ in $G$,
then $G$ admits an $L$-coloring
if and only if 
$G$ admits an $\M_{L}$-coloring.
This implies $\chi_{\ell}(G) \leq \chi_{\text{DP}}(G)$.
Dvo\v{r}\'{a}k and Postle \cite{DP} showed that $\chi_{DP}(G) \leq 5$ if $G$ is a  planar graph, and $\chi_{DP}(G) \leq 3$ if $G$ is a planar graph with girth at least 5.
Also, Dvo\v{r}\'{a}k and Postle \cite{DP} observed that $\chi_{DP}(G) \leq k+1$ if $G$ is $k$-degenerate.

There are infinitely many simple graphs $G$ 
satisfying $\chi_{\ell} (G) < \chi_{\text{DP}}(G)$:
It is known that
$\chi(C_n) = \chi_{\ell}(C_n) = 2 < 3 = \chi_{\text{DP}}(C_n)$
for each even integer $n \geq 4$.
Furthermore,
the gap $\chi_{\text{DP}}(G) - \chi_{\ell} (G)$
can be arbitrary large.
For example,
Bernshteyn \cite{Bernshteyn} showed that 
for a simple graph $G$ with average degree $d$,
we have $\chi_{\text{DP}}(G) = \Omega(d / \log d)$,
while Alon \cite{Alon} proved 
that $\chi_{\ell}(G) = \Omega(\log d)$ and 
the bound is sharp.
See \cite{BK2} for more detailed results.
Recently, there are some works on DP-colorings;
see \cite{Bernshteyn, BK, BKZ, DP,  KO}.

\subsection{DP-coloring vs signed coloring}
We here point out that
a signed coloring of a signed graph $(G,\sigma)$
is a special case of a DP-coloring of $G$.
Let $L$ be the list assignment of $G$
with $L(u) = N_k$ for any vertex $u$ in $G$.
Then for an edge $uv$ in $G$,
let 
$$M_{L,uv} = 
\begin{cases}
\big\{(u,i)(v,i) : i \in N_k\big\} & \text{ if $\sigma(uv) = 1$,}\\
\big\{(u,i)(v,-i) : i \in N_k\big\} & \text{ if $\sigma(uv) = -1$.}
\end{cases}
$$
With this definition,
it is easy to see that
the signed graph $(G,\sigma)$ admits a signed $k$-coloring
if and only if 
the graph $G$ admits an $\M_{L}$-coloring. 
Furthermore,
we see the similar relation for signed list-coloring.
Thus we have the following property.

\begin{prop}
\label{DP_sign_prop}
If $G$ is  DP-$k$-colorable, then the signed graph $(G, \sigma)$ is signed $k$-choosable for any sign function $\sigma$.
\end{prop}

In this paper we prove the following theorem.

\begin{thm} \label{main-thm}
For each $k \in \{3, 4, 5, 6\}$, every planar graph without $C_k$ is DP-4-colorable.
\end{thm}

Note that Theorem \ref{main-thm} is an extension of Theorems \ref{list-noC4} and Theorem \ref{signed-noC4} (by Proposition \ref{DP_sign_prop}).


\section{Proof of Theorem \ref{main-thm}}
\label{main-proof-sec}

We first show the case when $k \in \{3,5,6\}$.
In this case,
we can easily prove Theorem \ref{main-thm}
by using some known results.

\begin{lemma}
\label{356-lemma}
For each $k \in \{3, 5, 6\}$, if $G$ is a planar graph without $C_k$, then $G$ is 3-degenerate. 
\end{lemma}
\proof
When $G$ has no $C_3$, then the girth of $G$ is at least 4.  Thus,
Euler formula directly proves that $G$ is $3$-degenerate.  Lih and Wang \cite{LW02} showed that every planar graph without $C_5$ is $3$-degenerate.  And it was showed in \cite{FJMS} that every planar graph without $C_6$ is also 3-degenerate.
\qed

\begin{cor} For each $k \in \{3, 5, 6\}$, if $G$ is a planar graph without $C_k$, then $G$ is DP-4-colorable.
\end{cor}

Thus remaining case is when $G$ is a planar graph without $C_4$. 
Such graphs are not necessarily $3$-degenerate
(e.g.~consider the line graph of dodecahedral graph),
but we can instead use 
the following Lemma appeared in \cite{Lam}.

\emph{An $F_5^3$-subgraph} $H$ of a graph $G$ is a subgraph
isomorphic to the graph consisting of a 5-cycle and a 3-cycle that share an edge
 and satisfying $d_G(v) = 4$ for all vertex $v$. 
That is, $V(F_5^3) = \{v_1, v_2, v_3, v_4, v_5, v_6 \}$ and $v_1, v_2, v_3, v_4, v_5, v_6$ form the cycle $C_6$ with a chord $v_2v_6$ and $d_G(v_i) = 4$ for all $1 \leq i \leq 6$.

\begin{lemma} \label{main-lemma} (Lemma 1 in \cite{Lam})  
If $G$ is a planar graph without $C_4$, then $G$ contains an $F_5^3$-subgraph.
\end{lemma}

Now we are ready to prove the main theorem.

\noindent {\bf (Proof of Theorem \ref{main-thm})} \\
Let $G$ be a minimal counterexample to Theorem \ref{main-thm}.  That is, $G$ is a planar graph without $C_4$ and $G$ does not admit a DP-4-coloring, but any proper subgraph of $G$ admits a DP-4-coloring.  


Let $L$ be a list assignment of $G$,
and let $\M_{L}$ be a matching assignment over $L$.
By Lemma \ref{main-lemma},
$G$ contains an $F_5^3$-subgraph $H$.
Let $G' := G - V(H)$ and $L'(v) = L(v)$ for $v \in V(G')$.
By the minimality of $G$,
$G'$ admits an $\M_{L'}$-coloring.  
Thus there is an independent set $I'$ in the $\M_{L'}$-cover with $|I'| = |V(G)| - |V(F_5^3)| = |V(G)| - 6$.  For $v \in V(H) = \{v_1, v_2, v_3, v_4, v_5, v_6\}$, we define
\[
L^*(v) = L(v) \setminus \bigcup_{uv \in E(G)} \big\{c' \in L(v): (u,c)(v,c') \in M_{L,uv}  \mbox{ and } (u, c) \in I' \big\}.
\]
Then since $|L(v_i)| \geq 4 = d(v_i)$ for all $i$, we have that $|L^*(v_2)| \geq 3, \ |L^*(v_6)| \geq 3$, and $|L^*(v_j)| \geq 2$ for $j \in \{1, 3, 4, 5\}$.
We denote by $\M_{L^*}$
the restriction of $\M_{L}$ into $H$ and $L^*$.

\begin{claim} \label{key-claim}
The $\M_{L^*}$-cover has an independent set $I^*$ with $|I^*| = 6 = |V(F_5^3)|$.
\end{claim}
\proof
Since $|L^*(v_2)| \geq 3$ and $|L^*(v_1)| \geq 2$, we can color $c \in L^*(v_2)$ such that $L^*(v_1) \setminus \{c' : (v_2, c) (v_1, c') \in M_{L^*}\}$ has at least two available colors.
By coloring greedily in order $v_3, v_4, v_5, v_6, v_1$, we can find an independent set $I^*$ with $|I^*| = 6$.  This completes the proof of Claim \ref{key-claim}.
\qed

By Claim \ref{key-claim}, the $\M_{L}$-cover has an independent set $I = I' \cup I^*$ with $|I| = |I'| + |I^*| = |V(G)|$, which is a contradiction for the choice of $G$.  This completes the proof of Theorem \ref{main-thm}.

\bigskip
\noindent {\bf Acknowledgment.} We thank Professor Alexandr Kostochka for helpful comments.


\end{document}